\title{Characteristically simple Beauville groups, I:\\
cartesian powers of alternating groups}

\author{Gareth A. Jones\\
School of Mathematics\\
University of Southampton\\
Southampton SO17  1BJ, U.K.\\
{\tt G.A.Jones@maths.soton.ac.uk}
}

\documentclass[12pt]{article}
\usepackage{color}
\usepackage[latin1]{inputenc}
\usepackage{a4wide}
\usepackage{latexsym}
\usepackage{amsmath}
\usepackage{amssymb}
\usepackage{amsfonts}
\newtheorem{theorem}{Theorem}[section]
\newtheorem{lemma}[theorem]{Lemma}
\newtheorem{corollary}[theorem]{Corollary}
\newtheorem{proposition}[theorem]{Proposition}
\newtheorem{conjecture}[theorem]{Conjecture}
\date{}
\begin{document} 

\maketitle

\begin{abstract}
A Beauville surface (of unmixed type) is a complex algebraic surface which is the quotient of the product of two curves of genus at least 2 by a finite group $G$ acting freely on the product, where $G$ preserves the two curves and their quotients by $G$ are isomorphic to the projective line, ramified over three points. Such a group $G$ is called a Beauville group. We show that if a characteristically simple group $G$ is a cartesian power of a finite simple alternating group, then $G$ is a Beauville group if and only if it has two generators and is not isomorphic to $A_5$.
\end{abstract}


\section{Introduction}

A Beauville surface $\mathcal S$ of unmixed type is a complex algebraic surface of the form $({\mathcal C}_1\times{\mathcal C}_2)/G$, where ${\mathcal C}_1$ and ${\mathcal C}_2$ are complex algebraic curves (compact Riemann surfaces) of genera $g_i\ge 2$, and $G$ is a finite group acting faithfully on the factors ${\mathcal C}_i$ and freely on their product, with ${\mathcal C}_i/G\cong {\mathbb P}^1({\mathbb C})$ for $i=1, 2$ and the induced covering $\beta_i:{\mathcal C}_i\to{\mathbb P}^1({\mathbb C})$ ramified over three points. (We will not consider Beauville surfaces of mixed type, where elements of $G$ transpose the factors ${\mathcal C}_i$.) These surfaces have interesting rigidity properties, and since the first examples, with ${\mathcal C}_1={\mathcal C}_2$ Fermat curves, were introduced by Beauville~\cite[p.~159]{Bea}, they have been intensively studied by geometers such as Bauer, Catanese and Grunewald~\cite{BCG05, BCG06, Cat}.

Group theorists have recently considered which groups $G$, called {\sl Beauville groups}, and in particular which non-abelian finite simple groups, arise in this way~\cite{BBF, Fai, FMP, FG, FGJ, FJ, GP, GLL, GM}; see~\cite{JonCFI} for a survey. One easily shows that the smallest such group, the alternating group $A_5$, is not a Beauville group.  Bauer, Catanese and Grunewald~\cite{BCG05} conjectured that all other non-abelian finite simple groups are Beauville groups. This was proved with finitely many possible exceptions by Garion, Larsen and Lubotzky~\cite{GLL}, proved completely by Guralnick and Malle~\cite{GM}, and extended to quasisimple groups (with $A_5$ and $\hat{A_5}\cong SL_2(5)$ as the only exceptions) by Fairbairn, Magaard and Parker~\cite{FMP}. Our aim here is to begin a further extension of these results, to characteristically simple groups.

A finite group $G$ is characteristically simple (that is, it has no characteristic subgroups other than $G$ and $1$) if and only if it is isomorphic to a cartesian power $H^k$ of s simple group $H$~\cite[Satz I.9.12]{Hup}. (Finiteness is essential here: for instance, the additive group of any division ring is characteristically simple, since the multiplicative group acts transitively on the non-zero elements.) The abelian characteristically simple finite groups are the elementary abelian $p$-groups $G\cong C_p^k$, where $p$ is prime. Catanese~\cite{Cat} has shown that an abelian group $G$ is a Beauville group if and only if $G\cong C_n^2$ for some $n$ coprime to $6$, so the abelian characteristically simple Beauville groups are the groups $C_p^2$ for primes $p\ge 5$. We therefore restrict our attention to cartesian powers $G=H^k$ of non-abelian finite simple groups $H$. Since simple Beauville groups have been dealt with, we may assume that $k\ge 2$.

The ramification condition on the coverings $\beta_i$ implies that a Beauville group is a quotient of two triangle groups $\Delta_i$, so in particular it is a $2$-generator group. For any finite group $H$ and any integer $n\ge 1$, there are only finitely many normal subgroups $N$ in the free group $F_n$ of rank $n$ with $F_n/N\cong H$, so there is an integer $c_n(H)\ge 0$ such that $H^k$ is an $n$-generator group if and only if $k\le c_n(H)$. Then $c_2(H)$ is an upper bound on the values of $k$ such that $H^k$ is a Beauville group.


Computing $c_n(H)$ can be difficult, even for $n=2$, but if $H$ is a non-abelian finite simple group then $c_n(H)$ is equal to the number $d_n(H)$ of orbits of ${\rm Aut}\,H$ on ordered $n$-tuples which generate $H$. In some small cases, methods introduced by P.~Hall~\cite{Hal} allow $d_n(H)$ to be computed by hand: thus $c_2(A_5)=d_2(A_5)=19$, for example. More recently, Dixon~\cite{Dix05} and Mar\'oti and Tamburini~\cite{MT} have given asymptotic estimates and inequalities for $d_2(A_n)$. 

Our main result concerns the alternating groups $A_n$, which are simple for $n\ge 5$.

\begin{theorem}\label{mainthm}
Let $G=H^k$, where $H=A_n$ for $n\ge 5$, and $k\ge 1$. Then the following are equivalent:
\begin{itemize}
\item $G$ is a Beauville group;
\item $G$ is a $2$-generator group not isomorphic to $A_5$;
\item $n=5$ and $2\le k\le d_2(H)\;(=19)$, or $n>5$ and $k\le d_2(H)$.
\end{itemize}
\end{theorem}

The analogue of Theorem~\ref{mainthm} is also true when $H$ is a simple group $L_2(q)$, $L_3(q)$, $U_3(q)$, a Suzuki group $Sz(2^e)$, a `small' Ree group $R(3^e)$, or one of the $26$ sporadic simple groups~\cite{Jon13}. This suggests the following conjecture:

\begin{conjecture}
Let $G$ be a non-abelian finite characteristically simple group. Then $G$ is a Beauville group if and only if $G$ is a $2$-generator group and $G\not\cong A_5$.
\end{conjecture}

In order to prove that any group $G$ is a Beauville group, one needs to realise it as a quotient of two triangle groups $\Delta_i$, or equivalently to produce two generating triples $(a_i, b_i, c_i)$ ($i=1, 2$) for $G$, images of the canonical generators of $\Delta_i$, which satisfy $a_ib_ic_i=1$. If $G=H^k$ then for $i=1, 2$ the members of such a triple are $k$-tuples $a_i=(a_{ij})$, $b_i=(b_{ij})$ and $c_i=(c_{ij})$ such that, for $j=1,\ldots, k$, the elements $a_{ij}, b_{ij}$ and $c_{ij}$ of $H$ form $k$ generating triples for $H$ which satisfy $a_{ij}b_{ij}c_{ij}=1$ and are mutually inequivalent under the action of ${\rm Aut}\,H$. Provided $k\le d_2(H)$ there are $k$-tuples of such triples in $H$. However, in order for $G$ to be a Beauville group we need two such $k$-tuples, one for each $i$, satisfying the following extra conditions.

The curve ${\mathcal C}_i$ has genus $g_i\ge 2$ provided the orders $l_i, m_i$ and $n_i$ of $a_i, b_i$ and $c_i$ satisfy $l_i^{-1}+m_i^{-1}+n_i^{-1}<1$, a condition which is automatic unless $H\cong A_5$ and $k=1$. However, it is generally harder to satisfy the other condition, that $G$ acts freely on ${\mathcal C}_1\times{\mathcal C}_2$. The elements of $G$ with fixed points on ${\mathcal C}_i$ are the conjugates of the powers of $a_i, b_i$ and $c_i$, so for $i=1, 2$ these two subsets $\Sigma_i$ of $G$ must intersect in only the identity element. This is equivalent to $\Sigma_1\cap\Sigma_2$ containing no elements of order $p$ for each prime $p$ dividing $|G|$, a condition which can be ensured by a careful choice of the generating triples $(a_{ij}, b_{ij}, c_{ij})$ for $H$, and of how they are are arranged to form $k$ pairs for $j=1,\ldots, k$. Ensuring this generally requires separate arguments for small $k$ and for large $k$. Also, as is typical when considering infinite families of finite groups, small groups sometimes require special treatment.

\section{Beauville surfaces and structures}

A finite group $G$ is a Beauville group (of unmixed type) if and only if it has generating triples $(a_i, b_i, c_i)$ for $i=1, 2$, of orders $l_i, m_i$ and $n_i$, such that
\begin{enumerate}
\item $a_ib_ic_i=1$ for each $i=1, 2$,
\item $l_i^{-1}+m_i^{-1}+n_i^{-1}<1$ for each $i=1, 2$, and
\item no non-identity power of $a_1, b_1$ or $c_1$ is conjugate in $G$ to a power of $a_2, b_2$ or $c_2$.
\end{enumerate}

\noindent (In fact (2) is implied by (1) and (3), but it is useful to discuss it here.) We will call such a pair of triples $(a_i, b_i, c_i)$ a {\sl Beauville structure\/} for $G$. Property~(1) is equivalent to $G$ being a smooth quotient $\Delta_i/K_i$ of a triangle group $\Delta_i=\Delta(l_i, m_i, n_i)$ by a normal surface subgroup $K_i$ uniformising ${\mathcal C}_i$, with $a_i$, $b_i$ and $c_i$ the local monodromy permutations for the covering ${\mathcal C}_i\to {\mathcal C}_i/G$ at the three ramification points; we call $l_i, m_i$ and $n_i$ the {\sl periods\/} of $\Delta_i$. We call $(a_i, b_i, c_i)$ {\sl hyperbolic\/} if it satisfies~(2); this is equivalent to ${\mathcal C}_i$ having genus at least $2$, so that $\Delta_i$ acts on the hyperbolic plane $\mathbb H$, with ${\mathcal C}_i\cong{\mathbb H}/K_i$ and ${\mathcal C}_i/G\cong {\mathbb H}/\Delta_i\cong {\mathbb P}^1({\mathbb C})$; note that if  $l_i^{-1}+m_i^{-1}+n_i^{-1}\ge 1$ then any quotient of $\Delta_i$ is either solvable or isomorphic to $A_5$, so most of the groups we shall consider automatically satisfy~(2). Property~(3), which is generally the most difficult to verify, is equivalent to $G$ acting freely on  ${\mathcal C}_1\times{\mathcal C}_2$, since the elements with fixed points in ${\mathcal C}_i$ are the conjugates of the powers of $a_i, b_i$ and $c_i$.

By a {\sl triple\/} in a group $G$ we will mean an ordered triple $(a, b, c)$ of elements of $G$ such that $abc=1$; it is a {\sl generating triple\/} if $a, b$ and $c$ (and hence any two of them) generate $G$, and it has {\sl type\/} $(l, m, n)$ if $a, b$ and $c$ have orders $l, m$ and $n$ (known as its {\sl periods}). Two triples are {\sl equivalent\/} if an automorphism of $G$ takes one to the other; thus equivalent triples have the same type. A triple of type $(l, m, n)$ can be converted into one whose type is any permutation of $(l, m, n)$ by permuting (and if necessary inverting) its elements. We say that a Beauville structure, in this notation, has {\sl type\/} $(l_1, m_1, n_1; l_2, m_2, n_2)$; by the preceding remark we can permute or transpose the two types $(l_i, m_i, n_i)$ while preserving properties~(1), (2) and (3).

To show that a pair such as $a_1$ and $a_2$ satisfy~(3) it is sufficient to verify that, for each prime $p$ dividing both $l_1$ and $l_2$, the element $a_1^{l_1/p}$ of order $p$ is not conjugate to $a_2^{kl_2/p}$ for $k=1, 2, \ldots, p-1$; in particular, if $l_1$ is prime it is sufficient to verify this for $a_1$.  Similar remarks apply to every other pair chosen from the two triples.

For $i=1, 2$ let $\Sigma_i$ be the set of elements of $G$ conjugate to powers of $a_i, b_i$ or $c_i$; this is a union of conjugacy classes of $G$, closed under taking powers.  Condition~(3) is satisfied if and only if $\Sigma_1\cap\Sigma_2=\{1\}$. For any prime $p$, let $\Sigma_i^{(p)}$ denote the set of elements of order $p$ in $\Sigma_i$. The preceding remark then gives the following:

\begin{lemma}\label{primeorder}
If generating triples $(a_i, b_i, c_i)$ ($i=1, 2$) for a group $G$ have types $(l_i, m_i, n_i)$,  then they satisfy~$(3)$ if and only if $\Sigma_1^{(p)}\cap \Sigma_2^{(p)}=\emptyset$ for each prime $p$ dividing $l_1m_1n_1$ and $l_2m_2n_2$. In particular, they satisfy~$(3)$ if $l_1m_1n_1$ and $l_2m_2n_2$ are mutually coprime. \hfill$\square$
\end{lemma}

\section{Generating cartesian powers}

Here we consider whether the $k$-th cartesian power $H^k$ of a finite group $H$ can be a Beauville group for various $k>1$. The observation that any Beauville group must be a $2$-generator group immediately imposes restrictions on $k$, as follows.

If $H$ is a finite group, let
\[b(H)=\max\{k\mid H^k\; \hbox{is a Beauville group}\}\]
and let
\[c_2(H)=\max\{k\mid H^k\; \hbox{is a $2$-generator group}\},\]
with $b(H)=0$ or $c_2(H)=0$ if $H^k$ is not a Beauville group or a $2$-generator group for any $k\ge 1$. Then clearly
\[b(H) \le c_2(H),\]
so any upper bound on $c_2(H)$ is also an upper bound on $b(H)$.

\begin{lemma}
Let $H$ be a finite group. If $H$ is not perfect then $c_2(H)\le 2$.
\end{lemma}
\noindent{\sl Proof.}
Since $H$ is not perfect it has $C_p$ as an epimorphic image for some prime $p$, so $H^k$ has $C_p^k$ as an epimorphic image. If $H^k$ is a $2$-generator group, so is $C_p^k$, giving $k\le 2$. \hfill$\square$

\medskip

Thus in looking for Beauville groups among cartesian powers $H^k\;(k\ge 2)$, it is sufficient to consider arbitrary powers of perfect groups $H$, and cartesian squares of  imperfect groups $H$. (In the context of cartesian powers of finite groups, Wiegold~\cite{Wie} has shown that there is a more general dichotomy between perfect and imperfect groups $H$, with the rank $d(H^k)$ of $H^k$ respectively having essentially logarithmic or arithmetic growth as $k\to\infty$. However Tyrer Jones~\cite{Tyr} has constructed finitely-generated infinite groups satisfying $H\cong H^2$, so $d(H^k)$ is constant.)

When $H$ is finite and perfect, arguments due to P.~Hall~\cite{Hal} bound the values of $k$ one needs to consider. For any finite group $H$, let $d_2(H)$ denote the number of normal subgroups $N$ of the free group $F_2=\langle X, Y\mid -\rangle$ of rank $2$ with $F_2/N\cong H$, and let $\phi_2(H)$ be the number of $2$-bases (ordered generating pairs) for $H$. Any $2$-base $(x,y)$ for $H$ determines an epimorphism $\theta: F_2\to H$, sending $X$ and $Y$ to $x$ and $y$, so it determines a normal subgroup $N=\ker\theta$ of $F_2$ with $F_2/N\cong H$; conversely, every such normal subgroup arises in this way, with two $2$-bases yielding  the same normal subgroup $N$ if and only if they are equivalent under an automorphism of $H$. Thus $d_2(H)$ is the number of orbits of ${\rm Aut}\,H$ on $2$-bases for $H$. Since only the identity automorphism can fix a $2$-base, this action is semiregular, so we have

\begin{lemma}
If $H$ is a finite group then
\[\phi_2(H)=d_2(H)|{\rm Aut}\,H|.\]
\end{lemma}

\begin{corollary}
If $H$ is a non-identity finite group then
\[c_2(H)\le d_2(H)=\frac{\phi_2(H)}{|{\rm Aut}\,H|}.\]
\end{corollary}

\noindent{\sl Proof.} If $H^k$ is a $2$-generator group, then it is a quotient of $F_2$, so there are at least $k$ normal subgroups of $F_2$ with quotient group $H$. This gives $c_2(H)\le d_2(H)$, and Lemma~3.2 completes the proof. \hfill$\square$

\medskip

For some groups we have $c_2(H)<d_2(H)$: for instance, if $p$ is prime then $c_2(C_p)=2$ whereas $d_2(C_p)=p+1$. However, if $H$ is a non-abelian finite simple group we have equality in Corollary~3.3. This uses the following well-known result:

\begin{lemma}
Let $H$ be a non-abelian finite simple group, and let a group $\Gamma$ have distinct normal subgroups $N_1,\ldots, N_k$ with $\Gamma/N_i\cong H$. Then
$\Gamma/\cap_{i=1}^kN_i\cong H^k$.
\hfill$\square$
\end{lemma}

Applying this with $\Gamma=F_2$, we have the following:

\begin{corollary}
If $H$ is a non-abelian finite simple group then
\[c_2(H)=d_2(H)=\frac{\phi_2(H)}{|{\rm Aut}\,H|}.\]
\end{corollary}

Since $b(H)\le c_2(H)$, and since simple groups have been dealt with, this means that in order to prove Theorem~\ref{mainthm} it is sufficient to show that if $H=A_n$ then $H^k$ is a Beauville group for each $k=2,\ldots,d_2(H)$. We will do this by repeated use of the following immediate consequence of Lemma~3.4:

\begin{corollary}
Let $H$ be a non-abelian finite simple group. Then $k$-tuples $a=(a_j)$, $b=(b_j)$ and $c=(c_j)$ form a generating triple for $G=H^k$ if and only if their components $(a_j, b_j, c_j)$ for $j=1,\ldots, k$ form $k$ mutually inequivalent generating triples for $H$. \hfill$\square$
\end{corollary}

Although we have concentrated on $2$-generator groups, most of the results discussed here have obvious extensions to $n$-generator groups for all $n\in{\mathbb N}$~\cite[\S1.6]{Hal}.

\section{Evaluating $\phi_2(H)$, $d_2(H)$ and $c_2(H)$}

In~\cite{Hal}, Hall gave a method which, among many other applications, gives a formula for $\phi_2(H)$ for any finite group $H$. From this one can in theory deduce the value of $d_2(H)$, though in practice one usually has to be content with approximations.

Since any pair of elements of $H$ generate a unique subgroup $K$ we have
\[|H|^2=\sum_{K\le H}\phi_2(K).\]
Applying M\"obius inversion in the lattice $\Lambda$ of subgroups of $H$, we therefore have
\[\phi_2(H)=\sum_{K\le H}\mu(K)|K|^2\]
where $\mu$ is the M\"obius function for $\Lambda$. This is defined recursively by the formula
\[\sum_{L\ge K}\mu(L)=\delta_{K,H}\]
for each $K\le H$, where $\delta$ denotes the Kronecker delta-function.
For instance, Hall used this to show that $d_2(A_5)=19$ and $d_2(A_6)=53$.

For any non-abelian finite simple group $H$, we have
\[c_2(H)=d_2(H)=\frac{\phi_2(H)}{|{\rm Aut}\,H|}<\frac{|H|^2}{|{\rm Inn}\,H|.|{\rm Out}\,H|}
=\frac{|H|}{|{\rm Out}\,H|}.\]
Results of Dixon~\cite{Dix}, of Kantor and Lubotzky~\cite{KL}, and of Liebeck and Shalev~\cite{LS} show that a randomly-chosen pair of elements generate $H$ with probability approaching $1$ as $|H|\to\infty$, so this upper bound is asymptotically sharp, that is,
 \[d_2(H)\sim\frac{|H|}{|{\rm Out}\,H|}\quad\hbox{as}\quad |H|\to\infty.\]
The values of $|H|$ and $|{\rm Out}\,H|$ for all non-abelian finite simple groups $H$ are given in~\cite{ATLAS}. They show that for each of the infinite families of such groups, $|{\rm Out}\,H|$ grows much more slowly than $|H|$, so that in fact $d_2(H)$ grows almost as quickly as $|H|$. For instance, $d_2(A_n)\sim n!/4$ as $n\to\infty$  (see~\cite{Dix05, MT} for more detailed results).

\section{Beauville structures in cartesian powers}

We saw in Corollary~3.6 how to form generating triples in a cartesian power $G=H^k$ of a non-abelian finite simple group $H$. In this section we consider sufficient conditions for two such triples to satisfy the hypothesis $\Sigma_1^{(p)}\cap\Sigma_2^{(p)}=\emptyset$ of Lemma~\ref{primeorder} for all primes $p$, so that condition~(3) of a Beauville structure is satisfied.

For any prime $p$, and any $g=(g_j)\in G=H^k$, define the $p$-{\it profile} of $g$ to be the $k$-tuple $P_p(g)=(e_j)$, where $p^{e_j}$ is the highest power of $p$ dividing the order $o(g_j)$ of $g_j$. Define the $p$-{\it summit} of $g$ to be the set $S_p(g)$ of $j\in{\mathbb N}_k:=\{1,\ldots, k\}$ for which $e_j$ attains its maximum value, provided this is not $0$, and define $S_p(g)=\emptyset$ if $o(g)$ is not divisible by $p$. If an element $g'=(g'_j)\in G$ of order $p$ is conjugate to a power of $g$, then it has coordinates $g'_j$ of order $p$ at all $j\in S_p(g)$, with $g'_j=1$ elsewhere; in other words, the {\sl support\/} ${\rm supp}(g'):=\{j\mid g'_j\ne 1\}$ of $g'$ is equal to $S_p(g)$. If $T=(a, b, c)$ is a triple in $G$, define the $p$-{\it summit} $S_p(T)$ to be the set $\{S_p(a), S_p(b), S_p(c)\}$ of subsets of ${\mathbb N}_k$. The first part of the following result is now obvious, and the second follows from Lemma~\ref{primeorder}:

\begin{lemma}\label{disjointsummits}
Let $H$ be a finite group. If two triples $(a_i, b_i, c_i)$ ($i=1, 2$) in $G=H^k$ have disjoint $p$-summits $S_p(T_i)$ for some prime $p$, then $\Sigma_1^{(p)}\cap\Sigma_2^{(p)}=\emptyset$. If this happens for each prime $p$ dividing any of their periods, then $\Sigma_1\cap\Sigma_2=\{1\}$. \hfill$\square$
\end{lemma}

One way of ensuring that certain elements of ${\mathbb N}_k$ are or are not in $S_p(g)$, without needing to know all the coordinates of $g$, is to make at least one coordinate of $g$ $p$-{\it full\/}, meaning that its order is divisible by the highest power of $p$ dividing the exponent of $H$. In this case, $S_p(g)$ is the set of all $j$ such that $g_j$ is $p$-full.

For any triple $T=(a, b, c)$ in $H$, let $\nu_p(T)$ be the number of $p$-full elements among $a, b$ and $c$. Define two triples $T_i$ ($i=1, 2$) to be $p$-{\sl distinguishing} if $\nu_p(T_1)\ne\nu_p(T_2)$, and to be {\sl strongly} $p$-{\sl distinguishing} if, in addition, whenever $\nu_p(T_i)=0$ then either $p^2$ does not divide $\exp(H)$ or $p$ does not divide any of the periods of $T_i$.

\begin{lemma}\label{stronglypdist}
Suppose that a non-abelian finite simple group $H$ has a set ${\mathbb T}=\{(T_{1,s}, T_{2,s})\mid s=1,\ldots, t\}$ of pairs $(T_{1,s}, T_{2,s})$ of generating triples for $H$ such that
\begin{enumerate}
\item for each prime $p$ dividing $|H|$ there is some $s=s(p)\in\{1,\ldots, t\}$ such that $T_{1,s}$ and $T_{2,s}$ are a strongly $p$-distinguishing pair;
\item for each $i=1, 2$ the $3t$ triples consisting of $T_{i,1},\ldots, T_{i,t}$ and their cyclic permutations are mutually inequivalent.
\end{enumerate}
Then $G:=H^k$ is a Beauville group for each $k=3t,\ldots, d_2(H)$.
\end{lemma}

\noindent{\sl Proof.} Let $T_{i,s}=(x_{i,s}, y_{i,s}, z_{i,s})$ for $i=1, 2$ and $s=1,\ldots., t$. Define elements $a_i, b_i$ and $c_i$ of $G$ by using $x_{i,s}, y_{i,s}, z_{i,s}$ or $y_{i,s}, z_{i,s}, x_{i,s}$ or $z_{i,s}, x_{i,s}, y_{i,s}$ respectively in their $j$-th coordinate positions where $j=3s-2$, $3s-1$ or $3s$ for $s=1,\ldots, t$, so that their $j$-th coordinates form a cyclic permutation of $T_{i,s}$. Thus
\[a_i=(x_{i,1}, y_{i,1}, z_{i,1}, \ldots, x_{i,t}, y_{i,t}, z_{i,t}, \ldots),\]
\[b_i=(y_{i,1}, z_{i,1}, x_{i,1}, \ldots, y_{i,t}, z_{i,t}, x_{i,t}, \ldots),\]
\[c_i=(z_{i,1}, x_{i,1}, y_{i,1}, \ldots, z_{i,t}, x_{i,t}, y_{i,t}, \ldots).\]
If $k>3t$ then in coordinate positions $j$ for $j=3t+1,\ldots, k$ we use $k-3t$ further generating triples for $H$ which are inequivalent to each other and to those already used for $j=1,\ldots, 3t$. Thus $a_i, b_i$ and $c_i$ have $k$ mutually inequivalent generating triples for $H$ in their coordinate positions, so they form a generating triple for $G$. Since $t\ge 1$ we have $k\ge 3$, so $G\not\cong A_5$ and this triple is hyperbolic

Now let $p$ be a prime dividing $|H|$, so there is some $s=s(p)\in\{1,\ldots, t\}$ such that $T_{1,s}$ and $T_{2,s}$ are a strongly $p$-distinguishing pair. If $g\in\Sigma_i^{(p)}$ for some $i$ then ${\rm supp}(g)=S_p(d_i)$ where $d_i=a_i, b_i$ or $c_i$, so
\[|{\rm supp}(g)\cap\{3s-2, 3s-1, 3s\}|=\nu_p(T_{i,s}).\]
Since $\nu_p(T_{1,s})\ne\nu_p(T_{2,s})$ it follows that $g$ cannot be a member of $\Sigma_i^{(p)}$ for both $i=1$ and $i=2$. Thus $\Sigma_1^{(p)}\cap\Sigma_2^{(p)}=\emptyset$, as required.
\hfill$\square$

\medskip

If we assume only that $T_{1,s}$ and $T_{2,s}$ are $p$-distinguishing, we have the following:

\begin{lemma}\label{pdist}
Suppose that a non-abelian finite simple group $H$ has a set ${\mathbb T}=\{(T_{1,s}, T_{2,s})\mid s=1,\ldots, t\}$ of pairs $(T_{1,s}, T_{2,s})$ of generating triples for $H$ such that
\begin{enumerate}
\item for each prime $p$ dividing $|H|$ there is some $s=s(p)\in\{1,\ldots, t\}$ such that $T_{1,s}$ and $T_{2,s}$ are a $p$-distinguishing pair;
\item the $6t$ triples consisting of the $2t$ triples $T_{i,s}$ and their cyclic permutations are mutually inequivalent.
\end{enumerate}
Then $G:=H^k$ is a Beauville group for each $k=6t,\ldots, d_2(H)$.
\end{lemma}

\noindent{\sl Proof.} The proof is similar, with the first $3t$ coordinates of $a_i, b_i$ and $c_i$ defined as before. For $j=3t+1,\ldots, 6t$ we use the coordinate of $a_{3-i}, b_{3-i}$ or $c_{3-i}$ in position $j-3t$, and if $k>6t$ we  use further mutually inequivalent generating triples for $H$ in the remaining coordinate positions, so each triple $a_i, b_i, c_i$ generates $G$. For each prime $p$ dividing $|H|$, each of the six generators $a_i, b_i, c_i$ has at least one $p$-full coordinate, so again if $g\in\Sigma_i^{(p)}$ for some $i$ then ${\rm supp}(g)=S_p(d_i)$ where $d_i=a_i, b_i$ or $c_i$. The proof now continues as before. \hfill$\square$

\medskip

In many applications of these lemmas one can take $t=1$, so that a single pair of generating triples suffices. This often happens when $|H|$ is divisible by a small number of primes, as with the smaller alternating groups. The small values of $k$ not covered by these lemmas can often be dealt with by applying the following result:

\begin{lemma}\label{strongcoprimeperiods}
Let $H$ be a non-abelian finite simple group with $r$ mutually inequivalent generating triples of type $(l,m,n)$, where $r\ge 2$. If $l, m$ and $n$ are mutually coprime then $H^k$ is a Beauville group for each $k=2,\ldots, 6r$.
\end{lemma}

\noindent{\sl Proof.} Let the specified generating triples for $H$ be $(x_j, y_j, z_j)$ for $j=1, \ldots, r$. The $6r$ generating triples formed by cyclically permuting the entries of each $(x_j, y_j, z_j)$ and each $(z_j^{-1},y_j^{=1},x_j^{-1})$ are then mutually inequivalent. Since $H$ is non-abelian and simple and $G\not\cong A_5$, it follows that for each $k=2,\ldots, 6r$ one can form a hyperbolic generating triple for $H^k$ by using any set of $k$ of these triples in the different coordinate positions. In particular, one can choose two such triples of the forms
\[a_1=(x_1, x_2, \ldots),\quad b_1=(y_1, y_2,\ldots),\quad c_1=(z_1, z_2,\ldots)\]
and
\[a_2=(x_1, y_2, \ldots),\quad b_2=(y_1, z_2,\ldots),\quad c_2=(z_1, x_2,\ldots),\]
where in both cases the dots represent arbitrary choices of $k-2$ generating triples for $H$ from the remaining $6r-2$. If $l, m$ and $n$ are mutually coprime then any prime $p$ dividing $lmn$ must divide exactly one of $l, m$ and $n$. Thus if $g\in\Sigma_1^{(p)}$ then ${\rm supp}(g)$ contains both $1$ and $2$, whereas if $g\in\Sigma_2^{(p)}$ it contains only one of them. Thus $\Sigma_1^{(p)}\cap\Sigma_2^{(p)}=\emptyset$.
\hfill$\square$

\section{Primitive permutation groups}

In preparation for dealing with the groups $H=A_n$, here we present some results which imply that certain triples generate the alternating group. 

\begin{proposition}\label{primAn}
Let $H$ be a primitive permutation group of degree $n$. Then $H\ge A_n$ if any of the following conditions is satisfied:
\begin{enumerate}
\item $H$ has a subgroup with an orbit of length $m$, where $1<m<n/2$, fixing the remaining $n-m$ points;
\item $H$ contains an $m$-cycle, where $m=2,\ldots, n-3$;
\item $H$ contains a double transposition, with $n\ge 9$.
\end{enumerate}
\end{proposition}

\noindent{\sl Proof.} (1) This is Margraff's extension~\cite{Mar, Mar95} of a theorem of Jordan, that if $1<m\le n/2$ then $H$ is $3$-transitive (see also~\cite[Theorems 13.4, 13.5]{Wie}).

\noindent (2) This is a recent extension, by the author~\cite{Jon12}, of a theorem of Jordan which deals with the case where $m$ is prime. (Although Wielandt, in~\cite[Theorem 13.9]{Wie}, refers to Jordan's paper~\cite{Jor73} for this, it is not explicitly stated there. However, it follows easily from Th\'eor\`eme~I of his paper~\cite{Jor71}, together with Th\'eor\`eme~I of~\cite{Jor73}.)

\noindent(3) A proof of this can be found on Peter Cameron's blog~\cite{Cam}. The natural representation of $AGL_3(2)$ shows that the lower bound on $n$ cannot be relaxed.
\hfill$\square$

\medskip

Proposition~\ref{primAn}(2) follows from a more general result~\cite{Jon12}, classifying the primitive groups containing an $m$-cycle for any $m$; this uses the classification of finite simple groups,  and results of M\"uller~\cite{Mue96} and Feit~\cite{Fei} for $m=n-1$ and $m=n$.

In some situations, the assumption of primitivity can be replaced with transitivity, which is easier to verify, by using the following:

\begin{lemma}\label{transprim}
Let $H$ be a transitive permutation group of degree $n$, containing an $m$-cycle. If $m$ is coprime to $n$ and $m>n/2$, then $H$ is primitive. In particular, if $m$ is prime and $m>n/2$ then $H$ is primitive.
\end{lemma}

\noindent{\sl Proof.} For the first assertion, if $H$ is imprimitive, then being transitive it has blocks of the same size $b$ for some proper divisor $b$ of $n$. If the cycle acts trivially on the blocks, its support is contained in a single block, so $m\le b\le n/2$, against our hypothesis. If the cycle acts non-trivially on the blocks, its support is a union of blocks, so $b$ divides $m$, again contradicting our hypothesis. 

The second assertion follows immediately from the first if $m<n$, and if $m=n$ then $H$, as a transitive group of prime degree, must be primitive.
\hfill$\square$

\begin{corollary}\label{coprimecycleAn}
Let $H$ be a transitive permutation group of degree $n$, containing an $m$-cycle. If $m$ is coprime to $n$ and $n/2<m<n-2$, then $H\ge A_n$ 
\end{corollary}

\noindent{\sl Proof.}
This follows immediately from Lemma~\ref{transprim} and Proposition~\ref{primAn}(2).
\hfill$\square$

\medskip

The following result, concerning elements with two or three cycles, is also useful:

\begin{lemma}\label{coprimecycles}
Let $H$ be a transitive permutation group of degree $n$. If $H$ has an element $h$ with either of the following cycle structures, then $H\ge A_n$:
\begin{enumerate}
\item cycle structure $c, d$ for coprime integers $c, d>1$;
\item cycle structure $1, c, d$ for coprime integers $c, d>1$ such that neither $1+c$ nor $1+d$ divides $n$.
\end{enumerate} 
\end{lemma}

\noindent{\sl Proof.} (1) We first show that $H$ is primitive, so suppose that it is imprimitive, with blocks of size $b$, a proper divisor of $n$. Let the cycles of $h$ be $C$ and $D$, of lengths $c$ and $d$. If $C$ contains a block then it is a disjoint union of blocks, so $b$ divides $c$, and hence divides $n-c=d$, contradicting the fact that $c$ and $d$ are coprime. The same argument applies to $D$, so every block $B$ meets both $C$ and $D$. It follows that $\langle h\rangle$ permutes the blocks transitively, so they meet $C$ in the same number $r=|B\cap C|$ of points, giving $c=rb$. Similarly, the blocks all satisfy $|B\cap D|=s$ for some $s$, so $d=sb$. Thus $b$ divides both $c$ and $d$, again a contradiction. Hence $H$ is primitive.

Without loss of generality, suppose that $c<d$, so $c<n/2$ since $c+d=n$. Since $c$ and $d$ are coprime, $h^d$ is a cycle of length $c$. Since $H$ is primitive and $1<c<n/2$, Proposition~\ref{primAn}(1) implies that $H\ge A_n$.

\noindent(2) If $H$ were imprimitive, the block containing the unique fixed point of $h$ would be a union of cycles of $h$. This is impossible since neither $1+c$ not $1+d$ divides $n$, so $H$ is primitive. As in (1), Proposition~\ref{primAn}(1) completes the proof.
\hfill$\square$

\medskip

In this result, $h\in A_n$ if and only if $c\equiv d$ mod~$(2)$. Provided $n\ge 8$, if $n=2m$ is even, or if $n=2m+1$ is odd, we obtain an element $h\in A_n$ satisfying hypothesis~(1) or (2) respectively by taking $c, d=m\pm 1$ or $m\pm 2$ as $m$ is even or odd. 

The following result (see~\cite[Lemma~2.3]{Eve}) is also useful in proving primitivity:

\begin{lemma}\label{primitive} Let $H = \langle h_1,\ldots, h_r\rangle$ be a transitive permutation group containing a cycle $h$ of prime  length. Suppose that for each $i=1,\ldots, r$ there is an element of ${\rm supp}(h)$ whose image under $h_i$ is also in ${\rm supp}(h)$. Then $H$ is primitive. \hfill$\square$
\end{lemma}

\section{Small alternating groups}\label{smallalt}

We will now prove Theorem~\ref{mainthm}. The result is already known for $k=1$ (see \S 1), so it is sufficient to show that if $H=A_n$ with $n\ge 5$ then $H^k$ is a Beauville group for each $k=2,\ldots, d_2(H)$. In this section we will consider the groups $H=A_n$ for $n=5,\ldots, 11$. The remaining alternating groups will be considered in the next section, using general methods which do not always apply when $n$ is small.

In this section and the next we will consider various triples $(x, y, z)$ in $A_n$; the reader may find it helpful to represent these as directed graphs on $n$ vertices, with arcs corresponding to the actions of two of the generators (usually $x$ and $y$).

\begin{proposition}\label{Ansmalln}
If $n=5,\ldots, 11$ then $A_n^k$ is a Beauville group for each $k=2, \ldots, d_2(A_n)$.
\end{proposition}

\noindent{\sl Proof.} We will deal with these seven groups individually.

\subsection{$H=A_5$}

It is well known that $A_5$ is not a Beauville group: for instance, any generating triple for this group must contain an element of order $5$ (otherwise it generates a solvable and hence proper subgroup); Sylow's Theorems imply that all elements of order $5$ are conjugate to powers of each other, so no pair of generating triples can satisfy condition~$(3)$. Nevertheless, we will prove that all other $2$-generator cartesian powers of $A_5$ are Beauville groups.

Hall~\cite{Hal} showed that $d_2(A_5)=19$, giving $19$ equivalence classes of generating triples for $A_5$. 
Note that $A_5$ has two conjugacy classes of elements of order $5$, transposed under conjugation by odd permutations, and also under squaring.

\begin{proposition}\label{A_5^k}
The group $A_5^k$ is a Beauville group for each $k=2,\ldots, d_2(A_5)$.
\end{proposition}

\noindent{\sl Proof.}  The existence of Beauville structures for $k=3, \ldots, 19$ may be deduced immediately from Lemma~\ref{stronglypdist} with $t=1$, by using generating triples $T_i=(x_i, y_i, z_i)$ ($i=1, 2$) of types $(2, 5, 5)$ and $(3, 3, 5)$, such as
\[x_1=(1, 2)(3, 4),\quad y_1=(1, 4, 2, 3, 5),\quad z_1=(1, 5, 4, 2, 3)\]
and
\[x_2=(1, 2, 3),\quad y_2=(3, 4, 5),\quad z_2=(1, 3, 5, 4, 2).\]
The only maximal subgroup of $A_5$ containing a subgroup $\langle z_i\rangle\cong C_5$ is its normaliser, isomorphic to $D_5$, and in neither case does this contain $y_i$, so each $T_i$ generates $A_5$. Since $A_5$ has exponent $30$, this pair of triples are strongly $p$-distinguishing for each of the relevant primes $p=2, 3$ and $5$, so $A_5^k$ is a Beauville group for each $k=3,\ldots, 19$.

If $k=2$ we cannot use Lemma~\ref{strongcoprimeperiods}, since $A_5$ has no pairs of generating triples satisfying its hypotheses. Instead we use the following more specific approach.

Let $(a_{1j}, b_{1j}, c_{1j})$ be generating triples for $A_5$ of types $(5, 5, 5)$ and $(3, 5, 5)$ respectively for $j=1,2$, with all five generators of order $5$ conjugate in $A_5$. For instance, we could take
\[a_{11}=(1,2,3,4,5),\quad b_{11}=(1,4,5,2,3),\quad c_{11}=(1,2,4,5,3),\]
and
\[a_{12}=(1,2,4),\quad b_{12}=(1,2,3,4,5),\quad c_{12}=(1,5,2,4,3).\]
These triples are inequivalent, so the elements $a_1=(a_{11}, a_{12})$, $b_1=(b_{11}, b_{12})$ and $c_1=(c_{11}, c_{12})$ form a generating triple of type $(15, 5, 5)$ for the group $G=A_5^2$.

Now let $a_2=(a_{21}, a_{22})$, $b_2=(b_{21}, b_{22})$ and $c_2=(c_{21}, c_{22})$ in $G$, where $(a_{2j}, b_{2j}, c_{2j})$ is a generating triple for $A_5$ of type $(3, 5, 5)$ or $(5, 5, 5)$ respectively for $j=1,2$, with $b_{21}$ and $c_{21}$ conjugate to each other in $A_5$, but not conjugate to $a_{22}, b_{22}$ or $c_{22}$. For instance, we could take $(a_{21}, b_{21}, c_{21})=(a_{12}, b_{12}, c_{12})$ and $(a_{22}, b_{22}, c_{22})=(a_{11}, b_{11}, c_{11})^g$ for some odd permutation $g\in S_5$. Then $(a_2, b_2, c_2)$ is a generating triple for $G$, also of type $(15, 5, 5)$.

In order to use Lemma~\ref{disjointsummits} to show that these triples $(a_i, b_i, c_i)$ satisfy condition~(3), it is sufficient to consider the primes $p=3$ and $5$. Elements of $\Sigma_1^{(3)}$ or $\Sigma_2^{(3)}$ have their first or second coordinates respectively equal to the identity, so they cannot be equal. Elements of $\Sigma_1^{(5)}$ have the form $(g_1, g_2)$ with $g_1$ conjugate to $g_2$ or with $g_2=1$, whereas elements of $\Sigma_2^{(5)}$ have $g_1=1$ or $g_1$ not conjugate to $g_2$, so again they cannot be equal. Thus $\Sigma_1\cap\Sigma_2=\{1\}$, so these two triples form a Beauville structure for $G$.
\hfill$\square$

\subsection{$H=A_6$}

When $k=2$ we can use Lemma~\ref{strongcoprimeperiods}. Up to conjugacy in $S_6$, the group $H=A_6$ has four equivalence classes of generating triples $(x, y, z)$ of type $(3,4,5)$, given by taking $z$ to be a $5$-cycle $(a, b, c, d, e)$ and $x=(a, f, c)$, $(a, f, e)$, $(a, b, f)(c, d, e)$ or $(a, f, d)(b, c, e)$. Each of these triples generates $H$ since no maximal subgroup contains elements of orders $3, 4$ and $5$ (see~\cite{ATLAS}). Now ${\rm Aut}\,A_6$ contains $S_6$ with index $2$, and acts semi-regularly on generating triples, so it has two orbits on triples of this type. It follows from Lemma~\ref{strongcoprimeperiods} that $A_6^k$ is a Beauville group for $k=2,\ldots, 12$. 

For $k=3,\ldots, d_2(A_6)=53$ (see~\cite{Hal}) we can apply Lemma~\ref{stronglypdist} with $t=1$, using triples of types $(3, 5, 5)$ and $(4, 4, 5)$ such as
\[x_1=(1, 2, 3),\quad y_1=(1, 3, 4, 5, 6), \quad z_1=(1, 6, 5, 4, 2)\]
and
\[x_2=(1, 2, 3, 4)(5, 6),\quad y_2=(1, 3)(2, 5, 4, 6), \quad z_2=(1, 2, 3, 4, 5).\]
Each triple generates a subgroup of $A_6$ which is doubly transitive and therefore primitive, and which contains a $3$-cycle ($x_1$ or $x_2z_2^{-1}$ respectively), so by Proposition~\ref{primAn}(2) it generates $A_6$. Since $A_6$ has exponent $60$, these two triples are strongly $p$-discriminating for each of the relevant primes $p=2, 3$ and $5$.

\subsection{$H=A_7$}

For small $k$ we can use Lemma~\ref{strongcoprimeperiods}, with triples
\[x=(1, 2, 3),\quad y=(3, 4, 5, 6, 7),\quad z=(1, 3, 7, 6, 5, 4, 2) \]
and
\[x=(1, 2, 3)(4, 5, 6),\quad y=(1, 6, 7, 3, 4),\quad z=(1, 6, 3, 7, 5, 4, 2)\]
of type $(3, 5, 7)$. No maximal subgroup of $A_7$ contains elements of orders $3, 5$ and $7$, so they are both generating triples. They are inequivalent since ${\rm Aut}\,A_7=S_7$ and the elements of order $3$ have different cycle structures, so $A_7^k$ is a Beauville group for each $k=2,\ldots, 12$.

For larger $k$ we can use Lemma~\ref{stronglypdist} with $t=1$. We choose one of the above triples of type $(3, 5, 7)$, together with a triple of type $(4, 7, 7)$, such as
\[x=(1, 2, 3, 4)(5, 6),\quad y=(1, 5, 2, 4, 6, 7, 3),\quad z=(1, 2, 6, 3, 7, 5, 4).\]
This generates a transitive group, which is primitive since its degree is prime,
and must therefore be $A_7$ since it contains the $3$-cycle $xy^3=(1, 7, 5)$. Since $A_7$ has exponent $2^2.3.5.7$, these two triples deal with the cases $k=3,\ldots, d_2(A_7)$.

\subsection{$H=A_8,\ldots, A_{11}$}

The arguments when $n=8, 9, 10$ and $11$ are similar to those used for $A_6$ and $A_7$: when $k=3,\ldots, d_2(A_n)$ we apply Lemma~\ref{stronglypdist} with $t=1$, and when $k=2$ we apply Lemma~\ref{strongcoprimeperiods}. In each case, we will simply state some triples which can be used. Verifying that these are generating triples is straightforward, using results from the preceding section. Pairs of triples of the same type can easily be shown to be inequivalent by representing them as edge-labelled directed graphs on $n$ vertices, with the arcs representing the actions of $x$ and $y$: inequivalence of triples corresponds to non-isomorphism of their corresponding graphs.

\smallskip

When $H=A_8$, of exponent $2^2.3.5.7$, a pair of generating triples
\[x=(1,2)(3, 4, 5, 6),\quad y=(1, 4, 3, 7, 8),\quad z=(1, 8, 7, 6, 5, 4, 2),\]
\[x=(1, 2, 3),\quad y=(3, 2, 4, 5, 6, 7, 8), \quad z=(1, 3, 8, 7, 6, 5, 4),\]
 of types $(4,5,7)$ and $(3,7,7)$ satisfy the hypotheses of Lemma~\ref{stronglypdist}, so this deals with $k=3,\ldots, d_2(A_8)$. For $k=2$ we can use Lemma~\ref{strongcoprimeperiods}, with an inequivalent generating triple of type $(4, 5, 7)$, such as
\[x=(1,2)(3, 4, 5, 6),\quad y=(1, 7, 4, 3,  8),\quad z=(1, 8, 6, 5, 4, 7, 2).\]
.
 
When $H=A_9$, of exponent $2^2.3^2.5.7$, we can use generating triples
\[x=(1, 2, 3, 4)(5, 6, 7, 8),\quad y=(1, 4, 5, 8, 9),\quad z=(1, 9, 7, 6, 5, 3, 2),\]
\[x=(1, 2, 3, 4, 5),\quad y=(1, 2, 5)(3, 6, 7, 8, 9),\quad z=(1, 4, 3, 9, 8, 7, 6, 2, 5)\]
 of types $(4,5,7)$ and $(5, 15, 9)$ for $k=3,\ldots, d_2(A_9)$.
For $k=2$ we can use an inequivalent generating triple of type $(4, 5, 7)$, such as
\[x=(1, 2, 3, 4)(5, 6, 7, 8),\quad y=(1, 4, 3, 5, 9),\quad z=(1, 9, 8, 7, 6, 5, 2).\]

When $H=A_{10}$, of exponent $2^3.3^2.5.7$, we can use generating triples
\[x=(1, \ldots, 8)(9, 10),\quad y=(1, 8, 7, 9, 10),\quad z=(1, 9, 6, 5, 4, 3, 2),\]
\[x=(1, \ldots, 9),\quad y=(1, 3, 5, 7, 10),\quad z=(1, 10, 6, 5, 2)(3, 9, 8, 7, 4)\] of types $(8,5,7)$ and $(9, 5, 5)$ for $k=3,\ldots, d_2(A_{10})$. For $k=2$ we can use an inequivalent generating triple of type $(8, 5, 7)$, such as
\[x=(1, \ldots, 8)(9, 10),\quad y=(1, 8, 7, 6, 9),\quad z=(1, 10, 9, 5, 4, 3, 2).\]

When $H=A_{11}$, of exponent $2^3.3^2.5.7.11$, we can use generating triples
\[x=(1, 2, 3, 4, 5),\quad y=(2, 11)(3, 6, 7, 8, 9, 10, 5, 4),\quad z=(1, 5, 10, 9, 8, 7, 6, 2, 4),\]
\[x=(1, 2, 3, 4, 5, 6, 7),\quad y=(1, 6, 7, 2, 4, 5, 8, 9, 10, 11, 3),\quad z=(1, 2, 6, 7, 5, 3, 11, 10, 9, 8, 4)\]
 of types $(5, 8, 9)$ and $(7, 11, 11)$ for $k=3,\ldots, d_2(A_{11})$. For $k=2$ we can use two inequivalent generating triples of type $(11, 3, 8)$ given by taking
\[x=(1, 2, 3, 4, 5, 6, 7, 8, 9, 10, 11)\quad{\rm and}\quad
y=(1, 4, 2)\quad{\rm or}\quad (1, 10, 2).\]

\smallskip

Thus, for each $n=5,\ldots, 11$ we have shown that $A_n^k$ is a Beauville group for $k=2,\ldots, d_2(A_n)$.  This completes the proof of Proposition~7.1. 
\hfill$\square$

\section{Larger alternating groups}

If we try to apply the preceding method to $A_n$ for larger $n$, then as $n$ increases we need to consider more primes dividing the group order, namely, the $\pi(n)\sim n/\log n$ primes $p\le n$. The irregular distribution of primes means that the {\it ad hoc\/} approach used for $n=5,\ldots, 11$ will not work in general, so we need a more systematic method of proof to deal with $A_n$ for $n\ge 12$ (an assumption which will be maintained throughout this section).

\begin{proposition}\label{Ansmallk}
If $n\ge 12$ then $A_n^k$ is a Beauville group for each $k=2, \ldots, (n-5)(n-6)(n-7)/4$. 
\end{proposition}

\noindent{\sl Proof.} We will use Lemma~\ref{strongcoprimeperiods}, which requires a number $r\ge 2$ of inequivalent generating triples of a type with mutually coprime periods. For odd $n$ we take
\[x=(1, 2, \ldots, n-4) \quad{\rm and}\quad y=(s, n-3)(t, n-2)(u, n-1)(v, n)\]
where $s, t, u$ and $v$ are distinct elements of $\{1, \ldots, n-4\}$, so that $z:=(xy)^{-1}$ is an $n$-cycle. Since $x$ and $y$ are even permutations they generate a subgroup $H=\langle x, y\rangle$ of $A_n$. Since it contains $z$, $H$ is transitive, and since $n-4$ is coprime to $n$ and greater than $n/2$, Corollary~\ref{coprimecycleAn} implies that $H=A_n$. Thus $(x, y, z)$ is a generating triple of type $(n-4, 2, n)$ for $A_n$, with mutually coprime periods. By representing these triples as directed graphs on $n$ vertices, with arcs labelled $x$ or $y$, we see that their equivalence classes (under ${\rm Aut}\,H=S_n$) correspond to orbits of the additive group ${\mathbb Z}_{n-4}$ on its $4$-element subsets $\{s, t, u, v\}$. Since $n-4$ is odd, this action is semiregular, so the number of equivalence clases is ${n-4\choose 4}/(n-4)=(n-5)(n-6)(n-7)/24$.

The argument is similar for even $n$. We take
\[x=(1, 2, \ldots, n-3) \quad{\rm and}\quad y=(s, s+1)(t, n-2)(u, n-1)(v, n)\]
where $s, s+1, t, u$ and $v$ are distinct elements of $\{1,\ldots,n-3\}$, so that $z$ is an $(n-1)$-cycle fixing $s$. In this case $H=\langle x, y\rangle$ is doubly transitive, since it is transitive and the stabiliser of $s$ contains $z$, so it is primitive; since $x\in H$, Proposition~\ref{primAn}(2) implies that $H=A_n$. We therefore have ${n-5\choose 3}=(n-5)(n-6)(n-7)/6$ equivalence classes of generating triples of type $(n-3, 2, n-1)$, with mutually coprime periods.

In either case the number $r$ of equivalence classes is at least $(n-5)(n-6)(n-7)/24$, so Lemma~\ref{strongcoprimeperiods} implies that $A_n^k$ is a Beauville group for $k=2,\ldots, 6r$.
\hfill$\square$

\medskip

To deal with larger $k$ we will use Lemma~\ref{pdist}. If $n$ is large there is no single pair of generating triples which are $p$-distinguishing for all primes $p\le n$, as there is for $n=5,\ldots,11$; instead we will find one such pair $(T_p, T'_p)$ for each $p$, with $T_p$ having at least one $p$-full element, while $T'_p$ has none, so that the pair is $p$-distinguishing.

It is easy to find such generating triples $T_p$ and $T'_p$ for each $p$. However, in order to apply Lemma~\ref{pdist} we also need the triples $T_p, T'_p$ and their cyclic permutations to be mutually inequivalent. We do this by constructing the triples $T_p$ and $T'_p$ in such a way that the corresponding prime $p$ can be recognised from their three elements.

\begin{lemma}\label{pfulltripleAn}
If $n\ge 12$ there is a set ${\mathbb T}_n$ of generating triples $T_p$ for $A_n$, one for each prime $p\le n$, such that
\begin{enumerate}
\item at least one of the three elements of $T_p$ is $p$-full, and
\item the triples $T_p$ ($p\le n$) and their cyclic permutations are mutually inequivalent.
\end{enumerate}
\end{lemma}

\noindent{\sl Proof.} First let $n$ be even, say $n=2m$. If $p$ is odd take $x$ to be a cycle of length $p^e$, where $p^e\le n<p^{e+1}$, so that $x$ is $p$-full. Take $y$ to consist of two cycles of coprime lengths $m\pm 1$ or $m\pm 2$ as $m$ is even or odd, with both cycles meeting the support of $x$. Then $H=\langle x, y\rangle$ is a transitive subgroup of $A_n$, and applying Lemma~\ref{coprimecycles}(1) to $y$ shows that $H=A_n$. Given $x$, we can choose such an element $y$ so that $z:=(xy)^{-1}$ has at least two non-trivial cycles, so $x$ is the only cycle in the generating triple $T_p=(x, y, z)$. Thus $p$ can be recognised from $\{x, y, z\}$ as the only prime dividing the order of this element, so for the odd primes $p\le n$ the triples $T_p$ and their cyclic permutations are mutually inequivalent. If $p=2$ take $x$ to have two cycles of lengths $2^e$ and $2$ where $2^e+2\le n<2^{e+1}+2$ (so $e\ge 2$), together with $n-2^e-2$ fixed points, so $x$ is $2$-full. Take $y$ as before, ensuring that some cycle of $x$ meets both cycles of $y$, and that $z$ again has at least two non-trivial cycles. The resulting triple $T_2$ is not equivalent to any cyclic permutation of a triple $T_p$ for $p$ odd, since none of its members is a cycle, nor to any cyclic permutation of itself, since $x$ and $y$ have even and odd orders respectively.

If $n=2m+1$ is odd, let $x$ and $y$ be as above, but with $y$ having an extra fixed point, and all three of its cycles meeting the (or a) cycle of $x$. It is easy to check that $y$ satisfies the conditions of Lemma~\ref{coprimecycles}(2), so $\langle x, y\rangle=A_n$ and the proof proceeds as before. 
\hfill$\square$

\medskip

Although the cycle structures of $x$ and $y$ are completely specified in the above construction, there is usually some freedom of choice for that of $z$.

\begin{lemma}\label{coprimetripleAn}
If $n\ge 12$ there is a set ${\mathbb T}'_n$ of generating triples $T'_p$ for $A_n$, one for each prime $p\le n$, such that
\begin{enumerate}
\item none of the three elements of $T'_p$ is $p$-full, and
\item the triples $T'_p$ ($p\le n$) and their cyclic permutations are mutually inequivalent.
\end{enumerate}
\end{lemma}

\noindent{\sl Proof.} We shall construct the triples $T'_p$ so that the corresponding prime $p$ can be recognised from each of them, as follows. Given a prime $p$ such that $7\le p\le n$, let $q$ be the largest prime less than $p$ (so $q\ge 5$), and let $x$ be the $q$-cycle $(1, 2, \ldots, q)\in A_n$. We will choose $y\in A_n$ so that $\langle x, y\rangle=A_n$, giving a generating triple $T'_p=(x,y,z)$ for $A_n$ where $z=(xy)^{-1}$. We will arrange that $y$ is not a cycle, and if $z$ is a cycle then it is at least as long as $x$; thus $p$ is the smallest prime greater than the length of a shortest cycle among $x, y$ and $z$. The primes $p=2, 3$ and $5$ will be dealt with later, using separate constructions.

Let $p\ge 7$, and first suppose that $n$ is even, so $p<n$. We take $y=(1, 2)c_3\ldots c_r$, where $c_3,\ldots, c_r$ are disjoint cycles, and each $c_i$ has length $m_i$, joining $m_i-1$ points from $\{q+1,\ldots, n\}$ to the point $i\in\{3, \ldots, q\}$, so $r\le q$. We take $m_3+\cdots+m_r=n-q+r-2$, so all such points are joined, and $z=(xy)^{-1}$ is a cycle of length $n-1$, fixing $1$; thus $z$ is even, and hence so is $y$, so the triple $T'_p=(x, y, z)$ lies in $A_n$. It generates a transitive subgroup $H$, which is primitive by Lemma~\ref{primitive} with $h=x$, so Proposition~\ref{primAn}(2) gives $H=A_n$ since $q\le p-2<n-2$. This triple satisfies~(1) if $m_3,\ldots, m_r$ and $n-1$ are coprime to $p$, so that $x, y$ and $z$ have orders coprime to $p$. If $n\not\equiv 1$ mod~$(p)$ then provided $n\not\equiv q-1$ mod~$(p)$ we can take $r=3$, with a single cycle $c_3$ of length $n-q+1$ coprime to $p$; if $n\equiv q-1$ mod~$(p)$ we can take $r=4$, with cycles $c_3$ and $c_4$ of lengths $m_3=2$ and $m_4=n-q$ ($\equiv -1$ mod~$(p)$). (It is here that we need $p>5$, so that $q>3$, otherwise we cannot take $r=4$.) If $n\equiv 1$ mod~$(p)$ we instead define $y=(2,q)c_3\ldots c_r$, fixing $1$, with the cycles $c_i$ as above; then $z$ transposes $1$ and $q$, and acts as a cycle of length $n-2$ on the remaining points, so again $T'_p$ is contained in $A_n$, and generates this group. In this case we can take $r=3$, with $c_3$ a cycle of length $m_3=n-q+1\equiv 2-q$ mod~$(p)$, so $m_3$ is coprime to $p$ since $2<q<p$.

Now suppose that $n$ is odd. We take $y=(1, q, 2)c_3\ldots c_r$, with the cycles $c_i$ as above, so $z$ fixes $1$ and $q$ and acts an a cycle of length $n-2$ on the remaining points. Again, $T'_p$ generates $A_n$ (if $q=p-2=n-2$ we can take $r=4$ with $c_3$ and $c_4$ transpositions, so $y^2$ is a $3$-cycle and hence Proposition~\ref{primAn}(2) applies). This triple has periods coprime to $p$, and hence satisfies~(1), provided $n\not\equiv 2$ mod~$(p)$: we can choose $c_3$, or $c_3$ and $c_4$, as when $n$ is even, if $n\not\equiv q-1$ mod~$(p)$ or $n\equiv q-1$ mod~$(p)$ respectively. If $n\equiv 2$ mod~$(p)$ we instead define $y=(1, 2, q)c_3\ldots c_r$, with the cycles $c_i$ as before, so that $z$ is now an $n$-cycle; the same argument then applies.

Note that in all of the above triples $T'_p$ ($p\ge 7$), $x$ is a $q$-cycle for a prime $q$ such that $5\le q\le n-2$; $y$ has a cycle contained in ${\rm supp}(x)$, of length $2$ or $3$ as $n$ is even or odd, it has one or two more non-trivial cycles, each meeting ${\rm supp}(x)$ in one point, and it fixes the rest of ${\rm supp}(x)$; finally $z$ is a cycle of length $n, n-1$ or $n-2$, or has cycle structure $n-2, 2$. In particular, $x$ has a single non-trivial cycle, whereas $y$ has at least two. Given such a triple, or any cyclic permutation of it, one can use this general information to identify which elements of the triple are $x, y$ and $z$, and hence to identify the prime $p$ from the length $q$ of the cycle $x$.

When $p=5$ we will take $x=(1,2,\ldots, q)$ again, but with $q=7$ rather than $3$. (Using the preceding construction, and taking $q=3$, would force $r=3$, whereas we would need to allow $r=4$ for some $n$.) We must therefore ensure that $T'_5$ is not equivalent to $T'_{11}$ or its cyclic permutations, since this triple (but no other $T'_p$) uses the same element $x$. If $n$ is even we define $y=(2,6)c_3\ldots c_r$, with $r=3$ or $4$ and $c_i$ defined as before, so that $z$ has cycle structure $3, n-3$. We can do this, with each $m_i$ and $n-3$ coprime to $5$, provided $n\not\equiv 3$ mod~$(5)$; if $n\equiv 3$ mod~$(5)$ we take $y=(2,5)c_3\ldots c_r$ instead, so that $z$ has cycle structure $4, n-4$. Now suppose that $n$ is odd. If $n\not\equiv 2$ mod~$(5)$ we define $y=(1,2)(6,7)c_3\ldots c_r$, so that $z$ is an $(n-2)$-cycle, and if $n\equiv 2$ mod~$(5)$ we define $y=(1,2)(5,7)c_3\ldots c_r$, so that $z$ has cycle structure $1,2,n-3$. In all cases, applying Lemma~\ref{primitive} and Proposition~\ref{primAn}(2) to the $7$-cycle $x$ shows that the resulting triple $T'_5$ generates $A_n$. Moreover, if $n$ is even then $z$ has a cycle of length $3$ or $4$, while if $n$ is odd then $y$ has two non-trivial cycles contained in ${\rm supp}(x)$; in either case this shows that $T'_5$ is not equivalent to any cyclic permutation of $T'_{11}$, and hence of any $T'_p$ for $p>5$.

For $p=2$ or $3$, since $n\ge 12$ we can satisfy~(1) by imposing the weaker condition that none of $x, y$ and $z$ has a cycle of length divisible by $p^2$, rather than the more restrictive $p$.

Let $p=3$. For odd $n\not\equiv 0$ mod~$(9)$ let
\[x=(1,2,\ldots,n),\quad y=(1,2,3),\quad z=(1,2,n,n-1,\ldots,3).\]
Then $H$ is primitive since the only $\langle x\rangle$-invariant proper equivalence relations are given by congruence mod~$(m)$ for some proper divisor $m$ of $n$, and $y$ preserves none of these. Thus $H=A_n$ by Proposition~\ref{primAn}(2), giving a generating triple $T'_3$ of type $(n,3,n)$. For odd $n\equiv 0$ mod~$(9)$ let
\[x=(1,2,\ldots,n-4)(n-3,n-2,n-1),\quad y=(1,3,2)(n-4,n-3,n),\]
so that $z$ is an $(n-2)$-cycle. Then $H$ is $A_n$ by Corollary~\ref{coprimecycleAn} since it contains the $(n-4)$-cycle $x^3$ with $n-4>n/2$ and $\gcd(n,n-4)=1$, so we have a generating triple of type $(3(n-4),3,n-2)$. Now let $n$ be even. If $n\not\equiv 1$ mod~$(9)$ let
\[x=(1,2,\ldots,n-1),\quad y=(1,n,2),\]
so that $z$ is an $(n-1)$-cycle, giving a generating triple of type $(n-1,3,n-1)$. For even $n\equiv 1$ mod~$(9)$ let
\[x=(1,2,\ldots,n-3),\quad y=(1,n,2)(3,5,4)(n-3,n-2,n-1),\]
\[z=(2,n,n-3,n-1,n-2,n-4,n-5,\ldots,5).\]
Then $H=A_n$ by Corollary~\ref{coprimecycleAn} since $x\in H$, so we have a generating triple of type $(n-3,3,n-3)$. In all these cases, $y$ has order $3$, whereas no triple $T'_p$ for $p>3$ contains such an element.

Now let $p=2$. If $n$ is odd let
\[x=(1,\ldots, n-2),\quad y=(n-3,n-1)(n-2,n).\]
Then $H$ is primitive by Lemma~\ref{transprim} since $x\in H$, so $H=A_n$ by Proposition~\ref{primAn}(3) since $y\in H$. This gives a generating triple $T'_2$ of type $(n-2,2, n)$. If $n$ is even let
\[x=(1,\ldots, n-1),\quad y=(1,2)(n-1,n),\]
so a similar argument shows that this is a generating triple of type $(n-1, 2, n-1)$. In both cases, the presence of an element $y$ of order $2$ distinguishes $T'_2$ from all other triples $T'_p$ and their cyclic permutations, so conclusion~(2) is satisfied.
\hfill$\square$

\begin{lemma}\label{inequivalent}
If $n\ge 12$ then no triple in ${\mathbb T}_n$ is equivalent to a cyclic permutation of a triple in ${\mathbb T}'_n$.
\end{lemma}

\noindent{\sl Proof.}
Each triple $(x, y, z)\in{\mathbb T}_n$ contains an element $y$ with cycle structure $c, d$ or $1, c, d$ as $n$ is even or odd, where $c, d\ge 5$. However, if $(x', y', z')\in{\mathbb T}'_n$ then no element $x', y'$ or $z'$ has such a cycle structure: by the construction of $T'_p$, each of $x'$ and $z'$ is a single cycle (possibly with fixed points), or has a cycle of length $2$ or $3$, while $y'$ always has a cycle of length $2$ or $3$. Thus $(x, y, z)$ cannot be equivalent to a cyclic permutation of $(x', y', z')$.
\hfill$\square$

\begin{proposition}\label{Anlargek} If $n\ge 12$ then $A_n^k$ is a Beauville group for each $k=6\pi(n),\ldots, d_2(A_n)$. 
\end{proposition}

\noindent{\sl Proof.} The preceding three lemmas show that the $\pi(n)$ ordered pairs $(T_p, T'_p)$ of generating triples for $A_n$, where $p\le n$, satisfy the hypotheses of Lemma~\ref{pdist}.
\hfill$\square$

Proposition~\ref{Ansmallk} deals with $A_n^k$ where $k=2,\ldots, (n-5)(n-6)(n-7)/4$, while Proposition~\ref{Anlargek} does this for $k=6\pi(n),\ldots, d_2(A_n)$. Easy estimates show that
\[6\pi(n)\le 3(n+1)\le(n-5)(n-6)(n-7)/4\]
for all $n\ge 12$, so these two results cover each $k=2,\ldots, d_2(A_n)$ for such $n$. Combining this with Proposition~\ref{Ansmalln}, which deals with $n=5,\ldots, 11$, we have proved:

\begin{theorem}\label{A_n^k}
If $n\ge 5$ then $A_n^k$ is a Beauville group for each $k=2,\ldots, d_2(A_n)$. \hfill$\square$
\end{theorem}

This completes the proof of Theorem~\ref{mainthm}.

\end{document}